\documentclass[a4paper,tbtags,draft]{amsart}
\usepackage[T1]{fontenc}
\usepackage{fix-cm}
\usepackage{fourier}
\usepackage[scaled]{beramono}
\usepackage[all]{xy}

\theoremstyle{definition}

\theoremstyle{remark}

\theoremstyle{plain}

\DeclareMathOperator{\Hom}{Hom}
\DeclareMathOperator{\id}{1}
\DeclareMathOperator{\inj}{\iota}
\DeclareMathOperator{\eval}{ev}

\DeclareMathOperator{\Transpose}{tr}

\newsavebox{\mysavebox}
\newdimen\myvskip
\newcommand{\directdisplay}[1][0pt]{%
  \par\ \par\sbox{\mysavebox}{\textit{Proof.}\ \ }%
  \myvskip\ht\mysavebox\advance\myvskip\dp\mysavebox%
  \advance\myvskip\baselineskip\advance\myvskip\parskip%
  \advance\myvskip 1.9\abovedisplayskip%
  \advance\myvskip #1\vskip -\myvskip}

\newdir{ >}{{}*!/-5pt/{>}}
\newdir{ <}{{}*!/-5pt/{<}}
\newdir{ V}{{}*!/-5pt/{\vee}}
\newdir{ A}{{}*!/-5pt/{\wedge}}

\begin{document}

%--------------------
\title{Cartesian Closed Categories are Distributive}
\date{\today}

\author{Marco Benini}
\address{Dipartimento di Scienze Teoriche e Applicate\\
  Universit\`{a} degli Studi dell'Insubria\\
  via Maz\-zini 5, I-21100 Varese (VA), ITALY}
\email{marco.benini@uninsubria.it}
\urladdr{http://marcobenini.wordpress.com/}

\thanks{The author has been supported by the project \emph{Correctness
    by Construction} (CORCON), EU 7\textsuperscript{th} framework
  programme, grant n.~PIRSES-GA-2013-612638, and by the project
  \emph{Abstract Mathematics for Actual Computation: Hilbert's Program
    in the 21\textsuperscript{st} Century}, funded by the John
  Templeton's Foundation. Opinions expressed here and the
  responsibility for them are the author's only, not necessarily
  shared by the funding institutions.\\
  The author has also to thank the kind hospitality of
  prof.~H.~Ishihara and Dr.~T.~Nemoto when he was a visiting
  researcher at JAIST, where the results have been developed.}



\begin{abstract}
  A folklore result in category theory is that a (weakly) Cartesian
  closed category with finite co-products is distributive. Usually,
  the proof of this small result is carried on using the fact that the
  exponential functor is right adjoint to the product functor. And,
  since functors having a right adjoint preserve co-limits, the result
  follows immediately.  But, when we try to explicitly construct the
  arrows, things become a bit more involved. In rare cases, it is
  pretty useful to have the exact arrows that make this isomorphism to
  hold.  The purpose of this note is to develop the explicit proof
  with all the involved arrows constructed in an explicit way.
\end{abstract}
%--------------------

\maketitle

A (weakly) Cartesian closed category is a category having finite
products and exponentiation\footnote{Some authors define Cartesian
  closed categories as those having finite limits and
  exponentiation. In fact, the result here applies to both the weak
  and the strong version of the term, so we will freely use the term
  `Cartesian' with no further adjectives.}.  The statement says that
any Cartesian closed category $\mathbb{C}$ with finite co-products is
\emph{distributive}, that is, for any objects $A,B,C$, $A \times (B +
C) \cong (A \times B) + (A \times C)$.

One arrow is easy to construct: consider the following
commutative diagram
\begin{align*}
  \xymatrix@R=3ex{
    & (A \times B) + (A \times C) & \\
    A \times B \ar[ur]^{\inj_1} \ar[r]^{\pi_1} \ar[d]_{\pi_2} & A & A
    \times C \ar[ul]_{\inj_2} \ar[l]_{\pi_1} \ar[d]^{\pi_2} \\
    B \ar[dr]_{\inj_1} & A \times (B + C) \ar[u]^{\pi_1}
    \ar[d]_{\pi_2} & C \ar[dl]^{\inj_2} \\
    & B + C & }
\end{align*}

It immediately follows that there are two universal arrows to the
product $A \times (B + C)$:
\begin{align*}
  \xymatrix{
    & (A \times B) + (A \times C) & \\
    A \times B \ar[ur]^{\inj_1} \ar[r]^{\pi_1} \ar[d]_{\pi_2}
    \ar@{.>}[dr]|{\id_A \times \inj_1} & A & A
    \times C \ar[ul]_{\inj_2} \ar[l]_{\pi_1} \ar[d]^{\pi_2}
    \ar@{.>}[dl]|{\id_A \times \inj_2} \\ 
    B \ar[dr]_{\inj_1} & A \times (B + C) \ar[u]^{\pi_1}
    \ar[d]_{\pi_2} & C \ar[dl]^{\inj_2} \\
    & B + C & }
\end{align*}

So, the newly constructed arrows form a co-product diagram, which
enables to construct a co-universal arrow as below
\begin{align}\label{diagram:1}
  \vcenter{\xymatrix{
      & (A \times B) + (A \times C) \ar@{.>}@/_1.2em/[dd]|(.3){[\id_a 
        \times \inj_1, \id_a \times \inj_2]}& \\
      A \times B \ar[ur]^{\inj_1} \ar[r]^{\pi_1} \ar[d]_{\pi_2}
      \ar[dr]|{\id_A \times \inj_1} & A & A
      \times C \ar[ul]_{\inj_2} \ar[l]_{\pi_1} \ar[d]^{\pi_2}
      \ar[dl]|{\id_A \times \inj_2} \\ 
      B \ar[dr]_{\inj_1} & A \times (B + C) \ar[u]^{\pi_1}
      \ar[d]_{\pi_2} & C \ar[dl]^{\inj_2} \\
      & B + C & } }
\end{align}

Hence, the arrow $[\id_a \times \inj_1, \id_a \times \inj_2]\colon (A
\times B) + (A \times C) \to A \times (B + C)$ is the first half of
the isomorphism we are seeking. It should be noted that the
construction of such an arrow does not require
exponentiation.\vspace{1.5ex} 

The other side require exponentiation. In the first place, we have to
prove that the product functor $(A \times \mathord{-})$ is left
adjoint to the exponentiation functor $(\mathord{-})^A$.

To do this, we prove that $A^{(B \times C)} \cong \left(A^B\right)^C$
by explicitly constructing the isomorphism. Consider the commutative
diagram
\begin{align*}
  \xymatrix@C=5em{
    A^{(B \times C)} \times B \times C \ar[r]^-{\eval} & A \\
    \left(A^B\right)^C \times B \times C \ar[r]_-{\eval \times \id_B}
    & A^B \times B \ar[u]_{\eval} }
\end{align*}

Focusing on the top and the right arrows in the diagram, we notice an
evident exponential object: we call $\alpha$ the obvious
exponential transpose
\begin{align}\label{diagram:2}
  \xymatrix@C=5em{
    A^{(B \times C)} \times B \times C \ar[r]^-{\eval}
    \ar@{.>}[dr]|{\alpha \times \id_B} & A \\
    \left(A^B\right)^C \times B \times C \ar[r]_-{\eval \times
      \id_B} & A^B \times B \ar[u]_{\eval} }
\end{align}

Thus, the following diagram constructs another exponential transpose
$\gamma$ 
\begin{align}\label{diagram:3}
  \vcenter{\xymatrix@C=5em{
    A^{(B \times C)} \times B \times C \ar[r]^-{\eval} & A \\
    \left(A^B\right)^C \times B \times C \ar[r]_-{\eval \times \id_B}
    \ar@{.>}[u]|{\gamma \times \id_{B \times C}}
    & A^B \times B \ar[u]_{\eval} } }
\end{align}

Finally, we can construct the last exponential transpose $\delta$ as
follows 
\begin{align}\label{diagram:4}
  \vcenter{\xymatrix@C=5em{
    A^{(B \times C)} \times C \ar[dr]^{\alpha} \ar@{.>}[d]|-{\delta
      \times \id_C} \\
    \left(A^B\right)^C \times C \ar[r]_-{\eval}
    & A^B } }
\end{align}

It remains to prove that $\gamma$ and $\delta$ are inverses to each
other and thus isomorphisms, so $A^{(B \times C)} \cong
\left(A^B\right)^C$. To this aim, let us consider the diagram:
\begin{align*}
  \xymatrix@C=8em{
    A^{B \times C} \times B \times C \ar[r]^{\eval_2}
    \ar@{.>}@<-.5ex>[dr]_{\alpha \times \id_B}
    \ar@{.>}@<-.5ex>[d]_{\delta \times \id_{B} \times \id_C} & A \\
    \left(A^B\right)^C \times B \times C \ar@{.>}@<-.5ex>[u]_{\gamma
      \times \id_{B \times C}} \ar[r]_{\eval_3 \times \id_B} & A^B
    \times B \ar[u]_{\eval_1} } 
\end{align*}
where the arrows are named as before. Consider the arrow $\theta$:
\begin{align*}
  \theta \equiv & \eval_1 \circ (\eval_3 \times \id_B) \circ (\delta
  \times \id_B \times \id_C) \circ (\gamma \times \id_{B \times C}) \\
  = & \eval_1 \circ ((\eval_3 \circ (\delta \times \id_C)) \times
  \id_B) \circ (\gamma \times \id_{B \times C}) \\ 
  = & \eval_1 \circ (\alpha \times \id_B) \circ (\gamma \times \id_{B
    \times C}) \\
  = & \eval_2 \circ (\gamma \times \id_{B \times C}) \\
  = & \eval_1 \circ (\eval_3 \circ \id_B)
\end{align*}
thus, by uniqueness of the exponential transpose, $\Transpose \theta =
\gamma$, that is, 
\begin{align*}
  \id_{\left(A^B\right)^C \times B \times C} =
  (\delta \times \id_B \times \id_C) \circ (\gamma \times \id_{B \times
    C}) = (\delta \circ \gamma) \times \id_{B \times C}\enspace,
\end{align*}
so $\delta \circ \gamma = \id_{\left(A^B\right)^C}$.

In the opposite direction,
\begin{align*}
  \tau \equiv & \eval_2 \circ (\gamma \times \id_{B \times C}) \circ
  (\delta \times \id_B \times \id_C) \\
  = & \eval_1 \circ (\eval_3 \times \id_B) \circ (\delta \times \id_B
  \times \id_C) \\ 
  = & \eval_1 \circ ((\eval_3 \circ (\delta \times \id_C)) \times
  \id_B) \\
  = & \eval_1 \circ (\alpha \times \id_B) \\
  = & \eval_2
\end{align*}
thus, by uniqueness of the exponential transpose, as before,
\begin{align*}
  \id_{A^{B \times C} \times B \times C} = (\gamma \times \id_{B
    \times C}) \circ (\delta \times \id_B \times \id_C) = (\gamma
  \circ \delta) \times \id_{B \times C}\enspace,
\end{align*}
so, $\gamma \circ \delta = \id_{A^{B \times C}}$.\vspace{1.5ex}

The functor
\begin{align*}
  \begin{array}{rrcl}
    (A \times \mathord{-})\colon& \mathbb{C} & \to & \mathbb{C} \\
    & B & \mapsto & A \times B \\
    & f\colon B \to C & \mapsto & \id_A \times f\colon A \times B \to
    A \times C 
  \end{array}
\end{align*}
is left adjoint to the functor
\begin{align*}
  \begin{array}{rrcl}
    (\mathord{-})^A\colon& \mathbb{C} & \to & \mathbb{C} \\
    & B & \mapsto & B^A \\
    & f\colon B \to C & \mapsto & \Transpose(f \circ \eval)\colon B^A
    \to C^A  
  \end{array}
\end{align*}
where $\Transpose(f \circ \eval)$ is the exponential transpose of $f
\circ \eval$:
\begin{align*}
  \xymatrix{
    B^A \times A \ar[r]^-{\eval} \ar@{.>}[d]_{\Transpose(f \circ \eval) 
    \times \id_A} & B \ar[d]^{f} \\
    C^A \times A \ar[r]_-{\eval} & C }
\end{align*}
as it is immediate to see from the usual exponentiation diagram:
\begin{align*}
  \xymatrix{
    A \times B^A \ar[r]^-{\eval} & B \\
    A \times C \ar@{.>}[u]^{\id_a \times \Transpose{g}} \ar[ur]_{g} } 
\end{align*}

In the $\Hom$-sets notation, $\Hom(A \times B, C) \cong \Hom(A,
C^B)$. But, in this respect, it is useful to explicitly construct the
isomorphisms between the $\Hom$-sets. Evidently, given an arrow
$g\colon A \times B \to C$, the associated arrow becomes $\Transpose
g$. 

In the other direction, given an arrow $f\colon A \to C^B$, we can
construct the diagram:
\begin{align*}
  \xymatrix{
    A \ar[r]^{f} \ar@{<->}[d]_{\cong_1} & C^B \\
    1 \times A \ar@{.>}[d]_{\Transpose(f \circ \cong_1) \times \id_A}
    & \\ 
    \left(C^B\right)^A \times A \ar[ruu]_{\eval} \\
    C^{(A \times B)} \times A \ar@{<->}[u]_{\cong_2 \times \id_A} }
\end{align*}

Thus, the arrow $\theta$ associated to $f$ is obtained by composition
as in the diagram:
\begin{align}\label{diagram:5}
  \vcenter{\xymatrix@C=9em{
    A \times B \ar@{<->}[r]^-{\cong_3} \ar@{.>}[rrd]_{\theta}& 
    1 \times A \times B \ar[r]^-{(\cong_2 \circ \Transpose(f \circ
      \cong_1)) \times \id_{A \times B}} &
    C^{(A \times B)} \times A \times B \ar[d]^{\eval}\\
    & & C } }
\end{align}

Of course, $\cong_1$ and $\cong_3$ come from the following diagram
\begin{align*}
  \xymatrix{
    1 \times A \ar[r]^{\pi_2} \ar[d]_{\pi_1} & A \\
    1 & A \ar[l]^{!} \ar[u]_{\id_A} \ar@{.>}[ul]|{\langle !,
      \id_A\rangle} }
\end{align*}
and the dotted arrow has $\pi_2$ as an inverse, as it is immediate to
show.

So, we can finally construct the arrow $A \times (B + C) \to (A \times
B) + (A \times C)$. This is done by explicitly redoing the proof of
the dual of Proposition~3.2.2 in~\cite{BOR94A} in the case of
interest. Precisely, we prove that the functor $(A \times
\mathord{-})$ preserves the $B + C$ co-product.

The $B+C$ co-product corresponds to the diagram
\begin{align*}
  \xymatrix{
    & B + C & \\
    B \ar[ur]^{\inj_1} & & C \ar[ul]_{\inj_2} }
\end{align*}

So, the same diagram transformed via the $(A \times \mathord{-})$
functor yields the co-cone
\begin{align*}
  \xymatrix{
    & A \times (B + C) & \\
    A \times B \ar[ur]^{\id_A \times \inj_1} & & A \times C
    \ar[ul]_{\id_A \times \inj_2} }
\end{align*}

Consider any co-cone
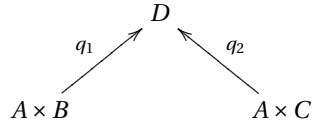
\begin{align*}
  \xymatrix{
    & D & \\
    A \times B \ar[ur]^{q_1} & & A \times C \ar[ul]_{q_2} }
\end{align*}

Thus, by exponentiation, i.e., the right adjoint of $(A \times
\mathord {-})$, 
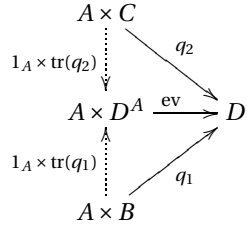
\begin{align*}
  \xymatrix{
    A \times C \ar[dr]^{q_2} \ar@{.>}[d]_{\id_A \times
      \Transpose(q_2)} & \\ 
    A \times D^A \ar[r]^{\eval} & D \\
    A \times B \ar[ur]_{q_1} \ar@{.>}[u]^{\id_A \times
      \Transpose(q_1)} & }
\end{align*}

So, the following is a co-cone
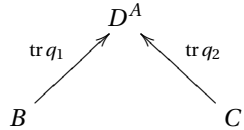
\begin{align*}
  \xymatrix{
    & D^A & \\
    B \ar[ur]^{\Transpose{q_1}} & & C \ar[ul]_{\Transpose{q_2}} } 
\end{align*}
and, being $B + C$ the co-product of $B$ and $C$, we get the
co-universal arrow $r\colon B + C \to D^A$
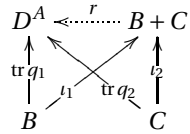
\begin{align*}
  \xymatrix{
    D^A & B + C \ar@{.>}[l]_{r} \\
    B \ar[u]|{\Transpose{q_1}} \ar[ur]|(.3){\inj_1} & C
    \ar[ul]|(.3){\Transpose{q_2}} \ar[u]|{\inj_2} }
\end{align*}

Applying the $(A \times \mathord{-})$ functor to this diagram, and
remembering that $\theta$ is the inverse of transpose, see
diagram~(\ref{diagram:5}), we get
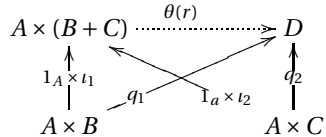
\begin{align*}
  \xymatrix@C=5em{
    A \times (B + C) \ar@{.>}[r]^{\theta(r)} & D \\
    A \times B \ar[u]|{\id_A \times \inj_1} \ar[ur]|(.3){q_1} & A
    \times C \ar[u]|{q_2} \ar[ul]|(.3){\id_a \times \inj_2} }
\end{align*}
which shows that $A \times (B + C)$ is the co-product of $A \times B$
and $A \times C$. Thus the arrow $A \times (B + C) \to (A \times B) +
(A \times C)$ is just the co-universal arrow of the $A \times (B + C)$
co-product, and it is evidently an inverse of the one synthesised in
diagram~(\ref{diagram:1}), because co-products are unique up to
isomorphisms. 

\bibliographystyle{amsalpha}
\bibliography{main}

\end{document}